\documentclass{amsart}
\DeclareMathSymbol{\twoheadrightarrow} 
{\mathrel}{AMSa}{"10}

\def\Q{{\mathbf Q}}
\def\Z{{\mathbf Z}}

\def\R{{\mathfrak  R}}
\def\F{{\mathbf F}}
\def\M{{\mathbf M}}
\def\L{{\mathbf L}}

\def\Sz{\mathbf{Sz}}

\def\Sn{{\mathbf S}_n}
\def\An{{\mathbf A}_n}
\def\bchi{{\mathbf \chi}}
\def\bphi{{\mathbf \phi}}
\def\Gal{\mathrm{Gal}}
\def\Perm{\mathrm{Perm}}
\def\tr{\mathrm{tr}}

\def\End{\mathrm{End}}
\def\Aut{\mathrm{Aut}}

\def\Mat{\mathrm{Mat}}
\def\cl{\mathrm{cl}}
\def\supp{\mathrm{supp}}

\def\I{\mathrm{Id}}

\def\fchar{\mathrm{char}}

\def\GL{\mathrm{GL}}

\def\SL{\mathrm{SL}}

\def\PSL{\mathrm{PSL}}

\def\M{\mathrm{M}}
\def\dim{\mathrm{dim}}

\def\P{{\mathbf P}}

\def\f{{\mathcal F}}

\newtheorem{thm}{Theorem}[section]
\newtheorem{lem}[thm]{Lemma}
\newtheorem{cor}[thm]{Corollary}

\theoremstyle{definition}
\newtheorem{defn}[thm]{Definition}
\newtheorem{ex}[thm]{Example}
\newtheorem{exs}[thm]{Examples}

\newtheorem{rem}[thm]{Remark}
\newtheorem{rems}[thm]{Remarks}
\title{Hyperelliptic jacobians and modular representations}
\author[Yuri G. Zarhin]{Yuri G. Zarhin}
\address{Department of Mathematics, Pennsylvania State University,
University Park, PA 16802, USA}
\email{zarhin\char`\@math.psu.edu}
\thanks{Partially supported by EPSRC grant GR/M 98135 and NSF grant DMS 0070664}

\begin{document}
\maketitle
\section{Introduction}

In \cite{Zarhin} the author  proved that
in characteristic $0$ the
jacobian $J(C)=J(C_f)$  of a hyperelliptic curve
$$C=C_f:y^2=f(x)$$
has only trivial endomorphisms
over an algebraic closure $K_a$ of the ground field $K$
if the Galois group $\Gal(f)$ of the irreducible polynomial
$f \in K[x]$ is ``very big". Namely, if $n=\deg(f) \ge 5$
and $\Gal(f)$ is either the symmetric group $\Sn$ or the alternating group $\An$
 then the ring $\End(J(C_f))$ of $K_a$-endomorphisms of $J(C_f)$ coincides with $\Z$. The proof was based on an explicit description of the Galois module $J(C_f)_2$
of points of order $2$ on $J(C_f)$. 
Namely, the action of the Galois group $\Gal(K)$ factors through $\Gal(f)$
 and the $\Gal(f)$-module $J(C_f)_2$ could be easily described in terms of the
 (transitive) action of $\Gal(f)$ on the set $\R_f$ of roots of $f$.

It turns out that if $\Gal(f)$ contains $\An$ then the Galois module $J(C_f)_2$ enjoys the following property (\cite{Zarhin}):

(*):each subalgebra in $\End_{\F_2}(J(C_f)_2)$ which contains the identity operator and is stable under the conjugation by Galois automorphisms either consists of scalars or coincides with
$\End_{\F_2}(J(C_f)_2)$. 

Applying (*) to the subalgebra $\End(J(C_f))\otimes \Z/2\Z$, one concludes that it consists of scalars, i.e.,
$\End(J(C_f)$ is a free abelian group of rank $1$
and therefore coincides with $\Z$. (The case
of $\End(J(C))\otimes \Z/2\Z=\End_{\F_2}(J(C_f)_2)$
could not occur  in  characteristic zero.)

The proof of (*) was based on the well-known explicit description of $J(C_f)_2$ \cite{Mumford}, \cite{MumfordENS} and elementary properties of $\An$ and its simplest nontrivial representation in characteristic $2$ of dimension $n-1$ or $n-2$ (depending on whether $n$ is odd or even).

In this paper we study  property (*) itself
from the point of view of representation theory over $\F_2$.
 Our results allow, in principle, to check the validity of (*)
 even if $\Gal(f)$ does not contain $\An$. 
We prove that $\End(J(C_f))=\Z$ for an infinite series 
of $\Gal(f)=\L_2(2^{r}):=\PSL_2(\F_{2^r})$ and $n=2^{r}+1$ 
(with $r \ge 3$ and $\dim(J(C_f))=2^{r-1}$) or when $\Gal(f)$ is 
 the Suzuki group $\Sz(2^{2r+1})$ and $n=2^{2(2r+1)}+1$ 
(with $\dim(J(C_f))=2^{4r+1}$).

We refer the reader to \cite{Mori1}, \cite{Mori2}, 
\cite{Katz1}, \cite{Katz2}, \cite{Masser}, \cite{Zarhin} 
for a discussion of known results about, and examples  
of, hyperelliptic jacobians without complex multiplication.

The paper is organized as follows. In \S \ref{mainr} 
 we state the main results and begin the discussion of 
 linear representations
for which an analogue of the property (*) holds true; 
 we call such representations {\sl very simple}. 
 In \S \ref{prdim3} we prove that the very simplicity of the 
 Galois module $X_{\ell}$ of points of prime order $\ell$ on 
 an abelian variety $X$ implies in characteristic zero that 
 $X$ does not have nontrivial endomorphisms.
In \S \ref{permute} we remind basic facts about permutation 
 groups and corresponding ordinary representations and modular 
 representations over $\F_2$. We use them in \S \ref{main2p} 
 in order to restate the main results as assertions about the 
 very simplicity of certain permutation modules using an 
 explicit description of points of order $\ell=2$ on 
 hyperelliptic jacobians. It turns out that  all 
 these permutations modules are Steinberg representations.
In \S \ref{St} we prove that the Steinberg representations 
are the only absolutely irreducible nontrivial representations 
(up to an isomorphism) over $\F_2$ for groups  $\L_2(2^r)$ 
 and $\Sz(2^{2r+1})$. 
In \S \ref{vsr} we  study very simple linear representations; 
in particular, we prove that all the (modular) Steinberg 
 representations discussed in Section \ref{St} are very simple. 
This ends the proof of main results.

\section{Main results}
\label{mainr}
Throughout this paper we assume that $K$ is a field. We fix its algebraic closure $K_a$
and write $\Gal(K)$ for the absolute Galois group $\Aut(K_a/K)$.
If $X$ is an abelian variety of dimension $g$ defined over $K$ then for each prime $\ell \ne \fchar(K)$ we write $X_{\ell}$ for the kernel of multiplication by $\ell$ in $X(K_a)$. It is well-known that $X_{\ell}$ is a $2g$-dimensional $\F_{\ell}$-vector space provided with a natural structure of $\Gal(K)$-module. We write $\End(X)$ for the ring of  $K_a$-endomorphisms of $X$ and $\End^0(X)$ for the corresponding finite-dimensional $\Q$-algebra $\End(X)\otimes\Q$.

The following notion plays a crucial role in this paper and will be discussed in detail in \S \ref{vsr}.

\begin{defn}
Let $V$ be a vector space over a field $\F$, let $G$ be a group and
$\rho: G \to \Aut_{\F}(V)$ a linear representation of $G$ in $V$. We
say that the $G$-module $V$ is {\sl very simple} if it enjoys the
following property:

If $R \subset \End_{\F}(V)$ be an $\F$-subalgebra containing the
identity operator $\I$ such that

 $$\rho(\sigma) R \rho(\sigma)^{-1} \subset R \quad \forall \sigma \in G$$
 then either $R=\F\cdot \I$ or $R=\End_{\F}(V)$.
\end{defn}

\begin{rems}
\label{image}
\begin{enumerate}
\item[(i)]
Clearly, the $G$-module $V$ is very simple if and only if the corresponding $\rho(G)$-module $V$ is very simple. 
\item[(ii)]
Clearly, if $V$ is very simple then the corresponding algebra homomorphism
	$$\F[G] \to \End_{\F}(V)$$
is surjective. Here $\F[G]$ stands for the group algebra of $G$. In particular, a 
very simple module is absolutely simple.
\item[(iii)]
If $G'$ is a subgroup of $G$ and the $G'$-module $V$ is very simple 
then the $G$-module $V$ is also very simple.
\item[(iv)]
Let $G'$ be a normal subgroup of $G$. If $V$ is a faithful 
very simple $G$-module then either $G' \subset \Aut_{\F}(V)$
consists of scalars (i.e., lies in $\F\cdot\I$) or the $G'$-module
$V$ is also very simple.
\end{enumerate}
\end{rems}

\begin{lem}
\label{dim3}
Let $X$ be an abelian variety of positive dimension $g$ over $K$. Let $\ell$ be a prime different from $\fchar(K)$. Assume that the $\Gal(K)$-module $X_{\ell}$ is very simple.
 Then either $\End(X)=\Z$ or $\fchar(K)>0$ and
$X$ is a supersingular abelian variety.
\end{lem}

We prove Lemma \ref{dim3} in \S \ref{prdim3}.

\begin{thm}
\label{maincSt}
Let $K$ be a field with $\fchar(K) \ne 2$,
$f(x) \in K[x]$ an irreducible separable polynomial 
 of degree $n\ge 5$. Let $\R=\R_f \subset K_a$ be 
 the set of roots of $f$, let $K(\R_f)=K(\R)$ be 
 the splitting field of $f$ and $\Gal(f):=\Gal(K(\R)/K)$ the Galois group of $f$, viewed as a subgroup of $\Perm(\R)$.
Let $C_f$ be the hyperelliptic curve $y^2=f(x)$. Let  $J(C_f)$ be
its jacobian, $\End(J(C_f))$ the ring of $K_a$-endomorphisms of $J(C_f)$.
Assume that $n$ and $\Gal(f)$ enjoy one of the following properties:

\begin{enumerate}
\item[(i)]
 $n=2^m+1 \ge 9$ and the Galois group $\Gal(f)$ of $f$ contains 
 a subgroup isomorphic to $\L_2(2^m)$;
\item[(ii)]
For some positive integer $k$ we have $n=2^{2(2k+1)}+1$ and the Galois group 
 $\Gal(f)$ of $f$ is isomorphic to $\Sz(2^{2k+1})$;
\end{enumerate}

Then:
\begin{enumerate}
\item[(a)]
The $\Gal(K)$-module $J(C)_2$ is very simple;
\item[(b)]
Either $\End(J(C_f))=\Z$ or 
$\fchar(K)>0$ and
$J(C_f)$ is a supersingular abelian variety.
\end{enumerate}
\end{thm}

\begin{rem}
It follows from  Lemma \ref{dim3} that in order to prove Theorem \ref{maincSt}, it suffices to check only the assertion a).
\end{rem}

\section{Proof of Lemma \ref{dim3}}
\label{prdim3}
Recall that 
$\dim_{\F_{\ell}}(X_{\ell})=2g$. 
Since $X$ is defined over $K$, one may associate with every $u
\in \End(X)$ and $ \sigma \in \Gal(K)$ an endomorphism
$^{\sigma}u \in \End(X)$ such that $$^{\sigma}u(x)=\sigma
u(\sigma^{-1}x) \quad \forall x \in X(F_a).$$
 Let us put
$$R:=\End(X) \otimes \Z/\ell\Z \subset \End_{\F_{\ell}}(X_{\ell}).$$
Clearly, $R$ satisfies all the conditions of Lemma \ref{dim3}.
This implies that either $R=\F_{\ell} \cdot \I$ or
 $R=\End_{\F_{\ell}}(X_{\ell})$.
If $\End(X) \otimes \Z/\ell\Z =R=\F_{\ell} \cdot \I$
 then the free abelian group $\End(X)$
has rank $1$ and therefore coincides with $\Z$. If $\End(X)
\otimes \Z/\ell\Z =R=\End_{\F_{\ell}}(X_{\ell})$ then the free abelian group
$\End(X)$ has rank $(2\dim(X))^2=(2g)^2$ and therefore
 the $\Q$-algebra $\End^0(X)$ has dimension $(2g)^2$.

Now Lemma \ref{dim3} becomes an immediate corollary of the following assertion proven in \cite{Zarhin}(see Lemma 3.1).

\begin{lem}
Let $Y$ be an abelian variety of dimension $g$ over an algebraically closed field $K_a$. Assume that  the semisimple $\Q$-algebra
$\End^0(Y)=\End(Y)\otimes\Q$ has dimension $(2g)^2$.
Then $\fchar(K_a)>0$ and $Y$ is supersingular.
\end{lem}

\section{Permutation groups and permutation modules}
\label{permute}

Let $B$ be a finite set consisting of $n \ge 5$ elements. We write $\Perm(B)$ for the group of permutations of $B$. A choice of ordering on $B$ gives rise to an isomorphism
$$\Perm(B) \cong \Sn.$$
Let $G$ be a  subgroup of $\Perm(B)$.
For each $b \in B$ we write $G_b$ for the stabilizer of $b$ in $G$; it is a subgroup of $G$.

\begin{rem}
\label{transitive}
Assume that the action of $G$ on $B$ is transitive.
It is well-known that each $G_b$ is a subgroup of index $n$ in $G$ and 
 all the $G_b$'s are conjugate one to another in $G$. 
 Each conjugate of $G_b$ in $G$ is the stabilizer of a point in $B$. 
 In addition, one may identify the $G$-set $B$ with the set of cosets $G/G_b$ 
with the standard action by $G$. 
\end{rem}

Let $\F$ be a field. We write $\F^B$
for the $n$-dimensional $\F$-vector space of maps $h:B \to \F$.
The space $\F^B$ is  provided with a natural action of $\Perm(B)$ defined
as follows. Each $s \in \Perm(B)$ sends a map
 $h:B\to \F$ into  $sh:b \mapsto h(s^{-1}(b))$. The permutation module $\F^B$
contains the $\Perm(B)$-stable hyperplane
$$(\F^B)^0=
\{h:B\to\F\mid\sum_{b\in B}h(b)=0\}$$
and the $\Perm(B)$-invariant line $\F \cdot 1_B$ where $1_B$ is the constant function $1$. The quotient $\F^B/(\F^B)^0$ is a trivial $1$-dimensional $\Perm(B)$-module.

Clearly, $(\F^B)^0$ contains $\F \cdot 1_B$ if and only if $\fchar(\F)$ divides $n$. If this is {\sl not} the case then there is a $\Perm(B)$-invariant splitting
$$\F^B=(\F^B)^0 \oplus \F \cdot 1_B.$$

Clearly, $\F^B$ and $(\F^B)^0$  carry natural structures of $G$-modules. Their characters depend only on characteristic of $\F$.

Let us consider the case of $\F=\Q$. Then
 the character of $\Q^B$ sends each $g \in G$ into the number of fixed points of $g$ 
 (\cite{SerreRep}, ex. 2.2, p. 12); it is called the {\sl permutation character}. Let us denote by 
$\bchi=\bchi_B:G \to \Q$
 the character of $(\Q^B)^0$.
It is known that the $\Q[G]$-module $(\Q^B)^0$ is absolutely simple 
 if and only if $G$ acts doubly transitively on $B$ (\cite{SerreRep}, ex. 2.6, p. 17).
 Clearly, $1+\bchi$ is the permutation character.

Now, let us consider the case of $\F=\F_2$.
It is well-known that one may view $\F_2^B$ as the $\F_2$-vector
 space of {\sl all} subsets of $B$ with symmetric difference as a sum.
Namely, a subset $T$ corresponds to its characteristic function
$\chi_T:B \to \{0,1\}=\F_2$ and  a function $h: B \to \F_2$
corresponds to its support $\supp(h)=\{x  \in B\mid
h(x)=1\}$. Under this identification each $s \in G \subset \Perm(B)$ sends
$T$ into $s(T)=\{s(b)\mid b \in T\}$.

Under this identification the hyperplane $(\F_2^B)^0$ corresponds to the $\F_2$-vector
 space of {\sl all} subsets of $B$ of {\sl even} cardinality with symmetric difference as a sum.

If  $n$ is even then let us define the $\Perm(B)$-module 
$$Q_B:=(\F_2^B)^0/(\F_2 \cdot 1_B).$$
If $n$ is odd then let us put
$$Q_B:=(\F_2^B)^0.$$
When $n$ is even, the quotient $Q_B$ corresponds to the $n-2$-dimensional $\F_2$-vector
 space of {\sl all} subsets of $B$ of {\sl even} cardinality with symmetric difference as a sum where each subset $T \subset B$ of even cardinality is identified with its complement $B\setminus T$.

\begin{rem}
Clearly,
$\dim_{\F_2}(Q_B)=n-1$ if $n$ is odd and $\dim_{\F_2}(Q_B)=n-2$ if $n$ is even.
In both cases $Q_B$ is a faithful $G$-module. 
\end{rem}

 Let $G^{(2)}$ be the set of $2$-regular elements of $G$. 
 Clearly, the Brauer character of the $G$-module $\F_{2}^B$  
 coincides with the restriction of $1+\bchi_B$ to $G^{(2)}$. 
 This implies easily that the Brauer character of the $G$-module $(\F_{2}^B)^0$ 
 coincides with the restriction of   $\bchi_B$ to $G^{(2)}$.

\begin{rem}
\label{Bcharacter}
 Let us denote by 
$\bphi_B=\bphi$
 the Brauer character of the $G$-module $Q_B$.
 One may easily check that $\bphi_B$ coincides with the restriction of 
 $\bchi_B$ to $G^{(2)}$ if  $n$ is odd and with the restriction of 
$\bchi_B-1$ to $G^{(2)}$ if  $n$ is even. 
\end{rem}

\begin{rem}
\label{oddeven}
Assume that $n=\#(B)$ is {\sl even}.
Let us choose $b \in B$ and let $G':=G_b$ and $B'=B\setminus\{b\}$. 
Then $n'=\#(B')=n-1$ is odd and there is a canonical isomorphism
 of $G'$-modules
$Q_{B'} \cong Q_B$
defined as follows.
First, there is a natural $G'$-equivariant embedding $\F_2^{B'} \subset \F_2^B$ which could be
obtained by extending each $h: B' \to \F_2$ to $B$ by letting $h(b)=0$.
Second, this embedding identifies $(\F_2^{B'})^0$ with a hyperplane of $(\F_2^{B})^0$ 
 which does not contain $1_B$. Now the desired isomorphism is given by the composition
$$Q_{B'}=(\F_2^{B'})^0 \subset (\F_2^{B})^0 \to
(\F_2^{B})^0/(\F_2 \cdot 1_B)=Q_B.$$
This implies that if the $G'$-module $Q_{B'}$ is very simple then the $G$-module $Q_B$ is also very simple.
\end{rem}

\begin{rem}
\label{St2}
Assume that $G$ acts on $B$ doubly transitively,
$\#(B)$ is odd and $\#(B)-1=\dim_{\Q}((\Q_B)^0)$ coincides with the largest power of $2$ dividing $\#(G)$. Then it follows from a theorem of Brauer-Nesbitt (\cite{SerreRep}, Sect. 16.4, pp. 136--137 ; \cite{Hump}, p. 249) that $Q_B$ is an absolutely simple $\F_2[G]$-module. In particular, $Q_B$ is  (the reduction of) the Steinberg representation \cite{Hump}.
\end{rem}

\section{Points of order $2$ on hyperelliptic jacobians}
\label{main2p}
We keep all notations of Section \ref{mainr}. In addition, we assume that $K$ is a field of characteristic different from $2$.
Let $C$ be a hyperelliptic curve over $K$ defined by an
equation $y^2=f(x)$ where $f(x) \in K[x]$ is a polynomial
of degree $n\ge 5$ without multiple roots. The rational function
$x \in K(C)$ defines a canonical double cover $\pi:C \to \P^1$.
Let $B'\subset C(K_a)$ be the set of ramification points
of $\pi$ (Weierstra{\ss } points). Clearly, the restriction of
$\pi$ to $B'$ is an injective map
$\pi:B' \hookrightarrow \P^1(K_a)$, whose image is either the 
 set $\R=\R_f$ of roots of $f$ 
if $n$ is even or the disjoint
union of $\infty$ and $\R$ if $n$ is odd. By abuse
of notation, we also denote by $\infty$ the ramification point
lying above $\infty$ if $n$ is odd and by $\infty_1$ and $\infty_2$ two unramified points lying above $\infty$ if $n$ is even. Clearly, if $n$ is odd then $\infty \in C(K)$. If $n$ is even
then the $2$-element set $\{\infty_1,\infty_2\}$ is stable under the action of $\Gal(K)$.

Let us put
  $$ B = \{ (\alpha,0) \mid f(\alpha) = 0\} \subset C(K_a).$$
Then $\pi$ defines a bijection between $B$ and
 $\R$ which
commutes with the action of  $\Gal(K)$. 
If $n$ is even then $B$ coincides with $B'$. 
 In the case of odd $n$ the set $B'$ is the disjoint union
 of $B$ and $\infty$.

\begin{thm}
\label{main2}
Suppose $n$ is an integer which is greater than or equal to $5$. Suppose $f(x) \in K[x]$ is a
separable polynomial of degree $n$, $\R \subset K_a$ the set of roots of $f$, let $K(\R)$ be the
splitting field of $f$ and $\Gal(f):=\Gal(K(\R)/K)$
the Galois group of $f$.

Suppose $C$ is the hyperelliptic curve $y^2=f(x)$ of genus $g=[\frac{n-1}{2}]$ over $K$. Suppose
 $J(C)$  is the jacobian of $C$ and $J(C)_2$ is
the group of its points of order $2$,  viewed as a
$2g$-dimensional $\F_2$-vector space provided with the natural
action of $\Gal(K)$.
Then the homomorphism $\Gal(K) \to \Aut_{\F_2}(J(C)_2)$ factors 
through the canonical surjection $\Gal(K) \twoheadrightarrow \Gal(K(\R)/K)=\Gal(f)$
and the $\Gal(f)$-modules $J(C)_2$ and $Q_{\R}$ are
isomorphic. In particular, the $G(K)$-module $J(C)_2$ is very simple if and only if the $\Gal(f)$-module $Q_B$ is very simple.
\end{thm}

\begin{rem}
\label{BR}
Clearly, $\Gal(K)$ acts on $B$ through the canonical surjective homomorphism
$\Gal(K)\twoheadrightarrow  \Gal(f)$,  because all points of $B$
are defined over $K(\R)$ and the natural homomorphism $\Gal(f) \to \Perm(B)$ is injective.
Clearly,  $\pi: B \to \R$ is a bijection of $\Gal(f)$-sets. This implies easily that the
$\Gal(f)$-modules $Q_B$ and $Q_{\R}$ are isomorphic. 
So, in order to prove Theorem \ref{main2} it suffices to check that the $\Gal(K)$-modules $Q_B$ and $J(C)_2$ are isomorphic.
\end{rem}

\begin{proof}[Proof of Theorem \ref{main2}]
  Here is a well-known explicit description of the
group $J(C)_2$ of points of order $2$ on  $J(C)$. Let us denote by $L$ the $K$-divisor $2(\infty)$ on $C$ if $n$ is odd and the
$K$-divisor $(\infty_1)+(\infty_2)$ if $n$ is even.
In both cases $L$ is an effective divisor of degree $2$.
Namely, let $T \subset B'$ be a subset of even cardinality. Then
(\cite{Mumford}, Ch. IIIa, Sect. 2, Lemma 2.4; \cite{MumfordENS},
pp. 190--191; see also \cite{Mori2}) the divisor $e_T=\sum_{P \in
T}(P) -\frac{\#(T)}{2}L$ on $C$ has degree $0$ and $2 e_T$ is
principal. If $T_1,T_2$ are two subsets of even cardinality in
$B'$ then the divisors $e_{T_1}$ and $e_{T_2}$ are linearly
equivalent if and only if either $T_1=T_2$ or $T_2=B'\setminus
T_1$. Also, if $T=T_1\triangle T_2$ then the divisor $e_T$ is
linearly equivalent to $e_{T_1}+e_{T_2}$.
 Hereafter we use the symbol $\triangle$ for the symmetric
 difference of two sets.
 Counting arguments imply easily that each point of $J(C)_2$ is the class of $e_T$ for some $T$. We know that such a choice is not unique. However, in the case of odd $n$ if we
demand that $T$ does not contain $\infty$ then such a choice
always exists and unique.
 This observation leads to a canonical
group isomorphism $$Q_B=(\F_2^B)^0 \cong J(C)_2, \quad T \mapsto \cl(e_T)$$
in the case of odd $n$. 
Here $\cl$ stands for the linear equivalence class of a divisor.
In the case of even $n$ we are still able to define
a canonical surjective
group homomorphism 
$$(\F_2^B)^0 \to J(C)_2, 
\quad T \mapsto \cl(e_T)$$
and one may easily check that the kernel of this
map is the line generated by the set $B$, i.e.,
the line generated by the constant function $1_B$.
This gives rise to the injective homomorphism
$$Q_B=(\F_2^B)^0/(\F_2 \cdot 1_B) \to J(C)_2,$$
which is an isomorphism, by counting arguments.
So, in both (odd and even) cases we get a canonical
isomorphism $Q_B \cong J(C)_2$, which obviously
commutes with the actions of $\Gal(K)$. In other
words, we constructed an isomorphism of $\Gal(K)$-modiles $Q_B$ and $J(C)_2$. In  light of Remark
 \ref{BR}, this ends the proof of Theorem \ref{main2}.
\end{proof}

Combining Theorem \ref{main2} and Lemma \ref{dim3} (for $\ell=2$),
we obtain the following corollary.
 
\begin{cor}
\label{mainvsc}
Let $K$ be a field with $\fchar(K) \ne 2$,
 $K_a$ its algebraic closure,
$f(x) \in K[x]$ an irreducible separable polynomial of degree $n\ge 5$. Let $\R=\R_f \subset K_a$ be the set of roots of $f$, let $K(\R_f)=K(\R)$ be the splitting field of $f$ and $\Gal(f):=\Gal(K(\R)/K)$ the Galois group of $f$, viewed as a subgroup of $\Perm(\R)$.
Let $C_f$ be the hyperelliptic curve $y^2=f(x)$. Let  $J(C_f)$ be
its jacobian, $\End(J(C_f))$ the ring of $K_a$-endomorphisms of $J(C_f)$.
Assume that the $\Gal(f)$-module $Q_{\R}$ is very simple. 
 Then either $\End(J(C_f))=\Z$ or 
$\fchar(K)>0$ and
$J(C_f)$ is a supersingular abelian variety.
\end{cor}

Notice that in order to prove Theorem \ref{maincSt}, it suffices to check the following statement.

\begin{thm}
\label{redp}
Let $n$ be a positive integer, $B$ a $n$-element set,
 $H \subset \Perm(B)$ a permutation group. Assume that $(n,H)$ enjoy one of the
following properties:

\begin{enumerate}
\item[(i)]
 $n=2^m+1 \ge 9$ and  $H$ contains 
 a subgroup isomorphic to $\L_2(2^m)$;
\item[(ii)]
For some positive integer $k$ we have $n=2^{2(2k+1)}+1$ and $H$
 contains a subgroup isomorphic to $\Sz(2^{2k+1})$;
\end{enumerate}
Then the $H$-module $Q_B$ is very simple.
\end{thm}

\begin{proof}[Proof of Theorem \ref{maincSt}  
 modulo Theorem \ref{redp}]
Let us put 
$$n=\deg(f), B=\R, H=\Gal(f).$$
 It follows from Theorem \ref{redp} that the 
 $\Gal(f)$-module $\Q_{\R}$ is very simple. 
 Now the result follows readily from Corollary \ref{mainvsc}.
\end{proof}

We prove Theorem \ref{redp} at the end of \S \ref{vsr}.

\section{Steinberg representation}
\label{St}
In this section we prove that the Steinberg representation is the only nontrivial absolutely irreducible representation over $\F_2$ (up to an isomorphism) of groups
$\L_2(2^m))$ and $\Sz(2^{2k+1})$.
We refer to \cite{Hump} for basic properties of Steinberg representations.

Let us fix an algebraic closure of $\F_2$ and denote it by $\f$. We write $\phi: \f \to \f$ for the Frobenius automorphism $x \mapsto x^2$. Let $q=2^m$ be a positive integral power of two. Then
the subfield of invariants of $\phi^m:\f \to \f$ is a finite field $\F_q$ consisting of $q$ elements.
Let  $q'$ be an integral positive power of $q$.
 If $d$ is a positive integer and $i$ is a non-negative integer then for each
matrix $u \in \GL_d(\f)$ we write $u^{(i)}$ for the matrix obtained by raising each entry of $u$ to the $2^i$th power.

Recall that an element $\alpha \in \F_q$ is called {\sl primitive} if $\alpha \ne 0$ 
and has multiplicative order $q-1$ in the cyclic multiplicative group $\F_q^*$.

\begin{lem}
\label{trace}
Let $q>2$, let $d$ be a positive integer and let $G$ be a subgroup of $\GL_d(\F_{q'})$.
Assume  that there exists an element
$u \in G\subset\GL_d(\F_{q'})$, whose trace $\alpha$ lies in $\F_q^*$  
 and has multiplicative order $q-1$.
Let $V_0=\f^d$ and let 
$\rho_0: G \subset \GL_d(\F_{q'}) \subset \GL_d(\f)=\Aut_{\f}(V_0)$
be the natural $d$-dimensional representation of $G$ over $\f$. For each positive integer $i<m$
we define a $d$-dimensional $\f$-representation 
$$\rho_i: G \to \Aut(V_i)$$
as the composition of 
$$G \hookrightarrow \GL_d(\F_{q'}), \quad x \mapsto x^{(i)}$$
 and the inclusion map 
$$\GL_d(\F_{q'}) \subset \GL_d(\f)\cong \Aut_{\f}(V_i).$$
Let $S$ be a subset of $\{0,1, \ldots m-1\}$.
Let us define a $d^{\#(S)}$-dimensional $\f$-representation $\rho_S$ of $G$ as the tensor product of representations $\rho_i$ for all $i \in S$. If $S$ is a proper subset of $\{0,1, \ldots m-1\}$  then there exists an element $u \in G$ such that the trace of $\rho_S(u)$ does not belong to $\F_2$. In particular, $\rho_S$ could not be obtained by extension of scalars to $\f$ from a representation of $G$ over $\F_2$.
\end{lem}

\begin{proof}
Clearly, 
$$\tr(\rho_i(u))=(\tr(\rho_0(u))^{2^i} \quad  \forall u \in G.$$
This implies easily that 
$$\tr(\rho_S(u))=\prod_{i\in S}\tr(\rho_i(u))= (\tr(\rho_0(u))^M$$
 where $M=\sum_{i \in S}2^i$. Since $S$ is a {\sl proper} subset of $\{0,1, \ldots m-1\}$, we have
$$0<M < \sum_{i=0}^{m-1} 2^i=2^{m}-1=\#(\F_q^*).$$
Recall that there exists $u \in G$ such that $\alpha=\tr(\rho_0(u))$ lies in $\F_q^*$
and the exact multiplicative order of $\alpha$ is
 $q-1=2^m-1$. 

This implies that $0 \ne \alpha^M \ne 1$. Since $\F_2=\{0,1\}$, we conclude that $\alpha^M \not\in \F_2$. Therefore
$\tr(\rho_S(u))=(\tr(\rho_0(u))^M=\alpha^M \not\in\F_2$.
\end{proof}

\begin{thm}
\label{L2}
Let $q \ge 8$ be a power of $2$ and $G=\L_2(q)=\PSL_2(\F_q)=\SL_2(\F_q)$.
Let $\rho: G \to \Aut(V)$ be an absolutely
irreducible faithful representation of $G$ over $\f$.
If the trace map $\tr_{\rho}: G \to \f$ takes on values in $\F_2$ then $\dim_{\f}(V)=q$. In particular, $\rho$ is the Steinberg representation of $G$.
\end{thm}

\begin{proof}
Let us put $q'=q$. We have
$$G=\SL_2(\F_q)\subset \GL_2(\F_q).$$
Clearly, for each $\alpha \in \F_q$ one may find
a $2\times 2$ matrix with determinant $1$ and trace $\alpha$. This implies that $G$ satisfies the conditions of Lemma \ref{trace}.

The construction described in Lemma \ref{trace} allows us to construct a $d^{\#(S)}$-dimensional $\f$-representation $\rho_S$ of $G$ for each subset
$S$ of of $\{0,1, \ldots m-1\}$. It is well-known (\cite{Brauer}, pp. 588-589) that $\rho_S$'s  exhaust the list of all absolutely irreducible $\f$-representations of $G=\SL_2(\F_q)$ 
 and therefore  $\rho$ is isomorphic to $\rho_S$ for some $S$. It follows from Lemma \ref{trace} that either $S$ is empty or $S=\{0,1, \ldots m-1\}$. The case of empty $S$ corresponds to the trivial $1$-dimensional representation. Therefore $S=\{0,1, \ldots m-1\}$ 
and $\rho$ is $2^m=q$-dimensional.
\end{proof}

Suppose $m=2k+1 \ge 3$ is an odd integer. Let $q=2^m=2^{2k+1}$ and $d=4$. Recall (\cite{HB}. pp. 182--194) that the {\sl Suzuki group} $\Sz(q)$ is the subgroup of $\GL_4(\F_q)$ generated by the matrices $S(a,b), M(\lambda), T$ defined as follows. For each $a, b \in \F_q$ the matrix $S(a,b)$ is defined by
\[S(a,b)=
	\begin{pmatrix}
1& 0& 0& 0\\
a& 1& 0& 0\\
b& a^{2^{k+1}}& 1& 0\\
a^{2^{k+1}+2}+ab+b^{2^{k+1}}& a^{2^{k+1}+1}+b& a& 1
\end{pmatrix}
\]
and for each $\lambda \in \F_q^*$ the matrix $M(\lambda)$ is defined by
\[M(\lambda)=
	\begin{pmatrix}
\lambda^{1+2^k}& 0& 0& 0\\
0& \lambda^{2^k}& 0& 0\\
0& 0& \lambda^{-2^k}& 0\\
0& 0& 0& \lambda^{1+2^k}
\end{pmatrix}.\]
The matrix $T$ is defined by
\[T=
	\begin{pmatrix}
	0& 0& 0& 1\\
	0& 0& 1& 0\\
	0& 1& 0& 0\\
	1& 0& 0& 0
\end{pmatrix}.\]

Notice that the trace of $S(0,b)T$ is $b^{2^{k+1}}$.
This implies easily that for each $\alpha \in \F_q$ one may find an element of $\Sz(q) \subset \GL_4(\F_q)$ with trace $\alpha$. This implies that
$$G=\Sz(q) \subset \GL_4(\F_q)$$
satisfies the conditions of Lemma \ref{trace}.
Notice also that $\#(\Sz(q))=(q^2+1)q^2(q-1)$ (\cite{HB}, p. 187).

\begin{thm}
\label{Sz} 
Let $\rho: \Sz(q) \to \Aut(V)$ be an absolutely
irreducible faithful representation of $\Sz(q)$ over $\f$.
If the trace map $\tr_{\rho}: \Sz(q) \to \f$ takes on values in $\F_2$ then $\dim_{\f}(V)=q^2$. In particular, $\rho$ is the Steinberg representation of $G$.
\end{thm}

\begin{proof}
Let us put $q'=q$. We know that
$G=\Sz(q)\subset \GL_4(\F_q)$ 
  satisfies the conditions of Lemma \ref{trace}.

The construction described in Lemma \ref{trace} allows us to construct a $4^{\#(S)}$-dimensional $\f$-representation $\rho_S$ of $G$ for each subset
$S$ of of $\{0,1, \ldots m-1\}$. It is known (\cite{Mart}, pp. 56--57) that $\rho_S$'s exhaust the list of all absolutely irreducible $\f$-representations of $G=\SL_2(\F_q)$ and therefore  $\rho$ is isomorphic to $\rho_S$ for some $S$. It follows from Lemma \ref{trace} that either $S$ is empty or $S=\{0,1, \ldots m-1\}$. The case of empty $S$ corresponds to trivial $1$-dimensional representation. Therefore $S=\{0,1, \ldots m-1\}$ 
and $\rho$ is $4^m=q^2$-dimensional.
\end{proof}

\begin{rem}
Assume that in the case \ref{redp}(i) (resp. \ref{redp}(ii))
that $H=\L_2(2^m))$ (resp. $\Sz(2^{2k+1})$). It follows from Remark \ref{St2} that $Q_B$ is the Steinberg representation of $H$.
\end{rem}

\section{Very simple representations}
\label{vsr}
\begin{exs}
\begin{enumerate}
\item[(i)]
If $\dim(V)=1$ then $V$ is always very simple.
\item[(ii)]
Assume that there exist $G$-modules $V_1$ and $V_2$ such that $\dim(V_1)>1, \dim(V_2)>1$ and the $G$-module $V$ is isomorphic to $V_1 \otimes_{\F} V_2$. Then $V$ is {\sl not} very simple. Indeed, the subalgebra
$$R=\End_{\F}(V_1)\otimes {\F}\cdot \I_{V_2}\subset
\End_{\F}(V_1)\otimes_{\F}\End_{\F}(V_2)=
\End_{\F}(V)$$
is stable under the conjugation by elements of $G$ but coincides neither with $\F\cdot \I$ nor with $\End_{\F}(V)$. (Here $\I_{V_2}$ stands for the identity operator in $V_2$.)

\item[(ii)bis]
Let $X \to G$ be a central extension of $G$.
Assume that there exist $X$-modules $V_1$ and $V_2$ such that $\dim(V_1)>1, \dim(V_2)>1$ and $V$, viewed as $X$-module, is isomorphic to $V_1 \otimes_{\F} V_2$. Then $V$ is {\sl not} very simple as an $X$-module. Since $X$ and $G$ have the same images in $\Aut_{\F}(V)$, the $G$-module $V$ is also not very simple. 
\item[(iii)]
Assume that there exists a subgroup $G'$ in $G$ of finite index $m>1$ and a $G'$-module $V'$ such that
the $\F[G]$-module $V$ is {\sl induced} by the $\F[G']$-module $V'$. (In particular, $m$ must divide $\dim(V)$.) Then $V$ is {\sl not} very simple. Indeed, one may view $W$ as a $G'$-submodule of $V$ such that $V$ coincides with
the direct sum $\oplus_{\sigma \in G/G'}\sigma W$. Let $R=\oplus_{\sigma\in G/G'}\End_{\F}(\sigma W)$ be the algebra  of all operators sending each $\sigma W$ into itself. Then $R$ 
is stable under the conjugation by elements of $G$ but  coincides neither with $\F\cdot \I$ nor with $\End_{\F}(V)$. 
\end{enumerate}
\end{exs}

\begin{ex}
Let $n \ge 5$ be an integer, $B$  a $n$-element
set. Suppose 
 $G$ is either $\Perm(B) \cong \Sn$ or the only subgroup in $\Perm(B)$ of index
  $2$ (isomorphic to $\An$). Then the $G$-module $Q_B$ is very simple. If $n$ is odd then this assertion is proven in \cite{Zarhin}, Th. 4.1.
If $n$ is even then $n \ge 6$, $n'=n-1 \ge 5$ is odd and the result follows from the odd case combined with Remarks \ref{image} and \ref{oddeven}.
\end{ex}

\begin{rems}
\label{factor0}
Assume that there exist $G$-modules $V_1$ and $V_2$ such that
 $\dim(V_1)>1, \dim(V_2)>1$
and the $G$-module $V$ 
is isomorphic to $V_1 \otimes_{\F} V_2$.

\begin{enumerate}
\item[(i)]
If $V$ is simple then both $V_1$ and $V_2$ are also simple. Indeed, if say, $V'$ is a proper $G$-stable subspace in $V_1$ then $V'\otimes_{\F}V_2$ is a proper $G$-stable subspace in $V_1\otimes_{\F}V_2=V$.
\item[(ii)]
If $V$ is  absolutely simple then both $V_1$ and $V_2$ are also absolutely simple. Indeed, assume that say, $R_1:=\End_{G}(V_1)$ has $\F$-dimension greater than $1$. Then
$\End_G(V)=\End_G(V_1\otimes_{\F}V_2)$ contains
$R_1 \otimes  \I_{V_2} \cong R_1$ and therefore also has dimension greater than $1$.
\end{enumerate}
\end{rems}

\begin{lem}
\label{induce}
Let $H$ be a group, $\F$ a field and $V$  a simple $\F[H]$-module of finite $\F$-dimension $N$.
Let $R \subset \End_{\F}(V)$ be an
$\F$-subalgebra containing the identity operator $\I$ and such that
 $$u R u^{-1} \subset R \quad \forall u \in H.$$
Then:
\begin{enumerate}
\item[(i)]
The faithful $R$-module $V$ is semisimple.
\item[(ii)]
Either the $R$-module $V$ is isotypic
or there exists a subgroup $H' \subset H$ of index $r$ dividing $N$ 
and a $H'$-module $V'$ of finite $\F$-dimension $N/r$ such that
	$r>1$ and
the $H$-module $V$ is induced by $V'$. In addition, if $\F=\F_2$ then $r<N$.
\end{enumerate}
\end{lem}

\begin{proof}

We may assume that $N>1$. Clearly, $V$ is a faithful $R$-module and
 $$u R u^{-1} = R \quad \forall u \in H.$$

{\bf Step 1}. $V$ is a {\sl semisimple} $R$-module. Indeed, let $U\subset V$ be a simple $R$-submodule. Then $U'=\sum_{s\in H} sU$
is a non-zero $H$-stable subspace in $V$ and therefore must
coincide with $V$. On the other hand, each $sU$ is also a
$R$-submodule in $V$, because $s^{-1} R s=R$. In addition, if $W\subset sU$ is an $R$-submodule then $s^{-1}W$ is an $R$-submodule
in $U$, because

$$R s^{-1}W=s^{-1} s R  s^{-1}W=s^{-1}RW=s^{-1}W.$$

 Since $U$ is simple, $s^{-1}W=\{0\}$ or $U$. This implies that $sU$ is also simple. Hence
 $V=U'$ is a sum of simple $R$-modules and therefore  is a  semisimple $R$-module.

{\bf Step 2}. The $R$-module $V$ is either isotypic or induced. Indeed, let us
split the semisimple $R$-module $V$ into the direct sum 
$$V =V_1\oplus \cdots \oplus V_r$$
 of its isotypic components.
 Dimension arguments imply that 
$r \le \dim(V) = N$.
 It follows easily from the arguments of the previous step  
that for each isotypic component $V_i$ its image $s V_i$ is an isotypic $R$-submodule 
 for each $s \in H$ and therefore is contained in some $V_j$.
  Similarly, $s^{-1}V_j$ is an isotypic submodule obviously containing $V_i$. 
 Since $V_i$ is the isotypic component,
  $s^{-1}V_j=V_i$ and therefore $sV_i=V_j$.
  This means that $s$ permutes the $V_i$; since $V$ is $H$-simple, 
$H$ permutes them transitively.

This implies that all $V_i$ have the same dimension $N/r$ and therefore 
 $r$ divides $\dim(V)=N$. Let $H'=H_i$ be the stabilizer of $V_i$ in $H$, i.e.
$$H_i =\{s \in H\mid sV_i=V_i\}.$$
The transitivity of the action of $H$ on $V_j$s implies that $[H:H_i]=r$.

If $r=1$ then $H=H'=H_i$.
This means that  $s V_i=V_i$ for all $s \in H$ and $V=V_i$ is isotypic.

Assume that $r>1$ and consider the $H'$-module $W=V_i$. 
 Clearly, $[H:H']=[H:H_i]=r$ divides $N$ and the $H$-module $V$ is iduced by $W$.

{\bf Step 3}. Assume that $r=N$ and $\F=\F_2$.
Then each $V_i$ is one-dimensional and contains exactly one non-zero vector say, $v_i$. Then the sum $\sum_{i=1}^N v_i$ is a non-zero $H$-invariant vector which contradicts 
 the simplicity of the $H$-module $V$.
\end{proof}

\begin{thm}
\label{central}
Suppose  $H$ is a group and 
$$\rho: H \to \Aut_{\F_2}(V)$$
 is an absolutely simple $\F_2[H]$-module of finite dimension $N$. Suppose there exists an $\F_2$-subalgebra $R \subset \End_{\F_2}(V)$ 
 containing the identity operator $\I$
 and such that
 $$u R u^{-1} \subset R \quad \forall u \in H.$$
Assume, in addition, that $H$ does not have nontrivial cyclic quotients of 
order dividing $N$.
If the $R$-module $V$ is isotypic then there exist 
 $\F_2[H]$-modules $V_1$ and $V_2$ such that $V$, viewed as $H$-module, is
 isomorphic to $V_1\otimes_{\F_2}V_2$ and the image 
of $R \subset \End_{\F_2}(V)$ under the induced isomorphism
$$\End_{\F_2}(V)=\End_{\F_2}(V_1\otimes_{\F_2} V_2)=
\End_{\F_2}(V_1)\otimes_{\F_2}\End_{\F_2}(V_2)$$
coincides with
$\End_{\F_2}(V_1)\otimes 
\I_{V_2}$. 
In particular, if both $V_1$ and $V_2$ have dimension greater than $1$ then the
 $H$-module $V$ is not very simple.
\end{thm}

\begin{proof}

 Since $V$ is isotypic, there exist a simple $R$-module $W$, a positive integer $d$ and an isomorphism
 $$\psi:V \cong W^d$$
of $R$-modules. Let us put
$$V_1=W, \quad V_2=\F_2^d.$$
The isomorphism $\psi$ gives rise to the isomorphism of $\F_2$-vector spaces
$$V=W^d=W\otimes_{\F_2}\F_2^d=
V_1\otimes_{\F_2} V_2.$$
We have
    $$d \cdot \dim(W)=\dim(V)=N.$$
Clearly, $\End_R(V)$ is isomorphic to the matrix algebra $\Mat_d(\End_R(W))$  of size $d$ over $\End_R(W)$.

Let us put
    $$k=\End_R(W).$$
Since $W$ is simple, $k$ is a finite-dimensional division algebra over $\F_2$. 
Therefore $k$ must be a finite field. 
 
We have

    $$\End_R(V) \cong \Mat_d(k).$$
Clearly, $[k:\F_2]$ divides $\dim_{\F_2}(W)$ and therefore divides $\dim_{\F_2}(V)=N$. 
Clearly, $\Aut(k/\F_2)$ is always a cyclic group of order $[k:\F_2]$ and therefore has order dividing $N$.

Clearly, $\End_R(V) \subset \End_{\F_2}(V)$ is stable under the adjoint action of $H$. This induces a homomorphism

 	$$\alpha: H \to\Aut_{\F_2}(\End_R(V))=\Aut_{\F_2}(\Mat_d(k)).$$

 Since $k$ is the center of $\Mat_d(k)$, it is stable under the action of $H$, i.e., we get a homomorphism $H\to\Aut(k/\F_2)$, which must be trivial, since $H$
is perfect and $\Aut(k/\F_2)$ is a cyclic group of order dividing $N$ and 
 therefore the kernel of the homomorphism  must coincide with $H$. This implies that the center $k$ of $\End_R(V)$ commutes with $H$. Since $\End_H(V)=\F_2$, we have $k=\F_2$. This implies that
$\End_R(V) \cong \Mat_d(\F_2)$ and
one may rewrite $\alpha$ as
 
$$\alpha: H \to \Aut_{\F_2}(\Mat_d(\F_2))=\Aut(\End_{\F_2}(V_2))=
\Aut_{\F_2}(V_2)/\F_2^*=\Aut_{\F_2}(V_2).$$

 It follows from the Jacobson density theorem that
 $R=\End_{\F_2}(W) \cong \Mat_m(\F_2)$ with $dm=N$. 

The adjoint action of $H$ on $R$ gives rise to a homomorphism
$$\beta: H \to \Aut_{\F_2}(\End_{\F_2}(W))=\Aut_{\F_2}(W)/\F_2^*=
\Aut_{\F_2}(W).$$
 
Clearly, $\alpha$ and $\beta$ provide $V_2$ and $V_1$ respectively
 with the structure of $H$-modules. Notice that
$$R=\End_{\F_2}(V_1)=\End_{\F_2}(V_1)\otimes
\I_{V_2}\subset\End_{\F_2}(V_1)\otimes_{\F_2}
\End_{\F_2}(V_2)=\End_{\F_2}(V).$$
Now our task boils down to comparison of the structures of $H$-module on
 $V=V_1\otimes_{\F_2}V_2$ defined by
$\rho$ and $\beta\otimes\alpha$ respectively.
I claim that
$$\rho(g)= \beta(g)\otimes \alpha(g) \quad \forall g \in H.$$
Indeed, 
notice that the conjugation by $\rho(g)$
in $\End_{\F_2}(V)=\End_{\F_2}(V_1\otimes_{\F_2}V_2)$ leaves stable
 $R=\End_{\F_2}(V_1)\otimes_{\F_2}\I_{V_2}$ and coincides
 on $R$ with the conjugation by $\alpha(g)\otimes\I_{V_2}$. Since the centralizer
of $\End_{\F_2}(V_1)\otimes
\I_{V_2}$ in  
$$\End_{\F_2}(V)=\End_{\F_2}(V_1)\otimes_{\F_2}
\End_{\F_2}(V_2)$$
coincides with 
$\I_{V_1}\otimes \End_{\F_2}(V_2)$,
there exists $u \in \Aut_{\F_2}(V_2)$ such that
$$\rho(g)=\beta(g)\otimes u.$$
Since the conjugation by $\rho(g)$ leaves stable
the centralizer of $R$, i.e. $\I_{V_1}\otimes \End_{\F_2}(V_2)$ 
and coincides on it with the conjugation by $\I_{V_1}\otimes \alpha(g)$, there exists a 
non-zero constant $\gamma \in \F_2^*$ such that
$u=\gamma \beta(g)$. This implies that
$$\rho(g)=\beta(g)\otimes u=\gamma\cdot \beta(g)\otimes\alpha(g).$$
Now one has only to recall that $\F_2^*=\{1\}$ and therefore $\gamma=1$.
\end{proof}

\begin{rem}
\label{factor2}
In the notations of Th. \ref{central} the $H$-modules $V_1$ and $V_2$ must be absolutely simple.
It follows easily from Remarks \ref{factor0}.
\end{rem}

Lemma \ref{induce} and Theorem \ref{central} together with Remark \ref{factor2} imply easily the following criterion of very simplicity over $\F_2$.

\begin{thm}
\label{very2}
Let $H$ be a group and $V$ be a $\F_2[H]$-module of finite dimension $N$ over $\F_2$.
Assume, in addition, that $H$ does not have nontrivial cyclic quotients of 
order dividing $N$ (e.g., $H$ is perfect).
 
Then $V$ is very simple if and only if the following conditions hold:
\begin{enumerate}
\item[(i)]
The $H$-module $V$ is absolutely simple;
\item[(ii)]
There do not exist a subgroup $H' \ne H$ of $H$ and a $\F_2[H']$-module $V'$ such that $V$ is induced by $V'$;
\item[(iii)]
There do not exist absolutely simple $\F_2[H]$-modules $V_1$ and $V_2$,
 both of dimension greater than $1$ and such that the $H$-module $V$ is isomorphic to $V_1\otimes_{\F_2} V_2$.
\end{enumerate}
\end{thm}

Combining Theorem \ref{very2} with Lemma \ref{induce} and Theorem \ref{central} we get easily the following corollary.

\begin{cor}
\label{Very2}
Let $H$ be a group and $V$ be a $\F_2[H]$-module of finite dimension $N$ over $\F_2$. Then $V$ is very simple if the following conditions hold:
\begin{enumerate}
\item[(i)]
The $H$-module $V$ is absolutely simple;
\item[(ii)]
$H$ does not contain a subgroup of finite index $r$
with $r\mid N$ and $1<r<N$.
 In addition, $H$ does not have cyclic quotients of order $N$,
 i.e., $H$ does not have a normal subgroup $H'$ of index $N$
with cyclic quotient $H/H'$;
\item[(iii)]
There do not exist absolutely simple $\F_2[H]$-modules $V_1$ and $V_2$, 
 both of dimension greater than $1$ and such that the $H$-module $V$ is isomorphic to $V_1\otimes_{\F_2} V_2$.
\end{enumerate}
\end{cor}

The following assertion follows easily from
Lemma \ref{induce} and Theorem \ref{central}.
 
\begin{cor}
\label{Very3}
Suppose  a positive integer $N>1$ and a group $H$ enjoy the following properties:
\begin{itemize}
\item
$H$ does not contain a subgroup of index dividing $N$ except $H$ itself.

\item
Let $N=ab$ be a factorization of $N$ into a product of two
positive integers $a>1$ and $b>1$. Then either
there does not exist an absolutely simple $\F_2[H]$-module of $\F_2$-dimension $a$ or
there does not exist an absolutely simple $\F_2[H]$-module of 
 $\F_2$-dimension $b$.
\end{itemize}

Then each absolutely simple $\F_2[H]$-module of $\F_2$-dimension $N$ is very simple. In other words, in dimension $N$ the properties of
absolute simplicity and very simplicity over $\F_2$ are equivalent.
\end{cor}

The next two theorems provide examples of very simple Steinberg representations.

\begin{thm}
\label{L2QB}
Let $q=2^m \ge 8$ be an integral power of $2$, let $B$ be a $(q+1)$-element set. Let $G'$ be a group acting faithfully on $B$. Assume that $G'$ contains a subgroup $G$ isomorphic to $\L_2(q)$. Then the $G'$-module $Q_B$ is very simple.
\end{thm}

\begin{proof}
We have
$\L_2(q)=G\subset G' \subset \Perm(B)$.
Clearly, it suffices to check that the $\L_2(q)$-module $Q_B$ is very simple.

First, notice that $\L_2(q)$ acts doubly transitively on $B$. Indeed, each subgroup 
of $\L_2(q)$ (except $\L_2(q)$ itself) has index $\ge q+1=\#(B)$ (\cite{Suzuki}, (6.27), p. 415). This implies that $\L_2(q)$ acts transitively on $B$. If the stabilizer $G_b$ of a point $b \in B$ has index $q+1$ then it follows easily from Th. 6.25 on p. 412 of \cite{Suzuki})
that $G_b$ in conjugate to the (Borel) subgroup of upper-triangular matrices and therefore the $\L_2(q)$-set $B$ is isomorphic to the projective line $\P^1(\F_q)$ with the standard action of $\L_2(q)$ which is well-known to be doubly (and even triply) transitive. By Remark \ref{St2}, this implies that the $\F_2[\L_2(q)]$-module $Q_B$ is absolutely simple. Recall that
$$\dim_{\F_2}(Q_B)=\#(B)-1=q=2^m.$$
By Theorem \ref{L2}, there no absolutely simple nontrivial $\F_2[\L_2(q)]$-modules 
 of dimension $< 2^m$. This implies that $Q_B$ is {\sl not} isomorphic to a tensor product of absolutely simple $\F_2[\L_2(q)]$-modules of dimension $>1$. Recall that all subgroups in $\L_2(q)$ different from $\L_2(q)$ have index $\ge q+1> q=\dim_{\F_2}(Q_B)$. It follows from Corollary \ref{Very2} that the $G$-module $Q_B$ is very simple. Since $G \subset G'$, the $G'$-module $Q_B$ is also very simple.
\end{proof}

\begin{thm}
\label{SzQB}
Let $k$ be a positive integer and $q=2^{2k+1}$, let $B$ be a $(q^2+1)$-element set. Let $G'$ be a group acting faithfully on $B$. Assume that $G'$ contains a subgroup $G$ isomorphic to $\Sz(q)$.  Then the $G'$-module $Q_B$ is very simple.
\end{thm}

\begin{proof}
We have
$\Sz(q) =G \subset G'\subset\Perm(B)$.
First, notice that $\Sz(q)$ acts doubly transitively on $B$. 
Indeed, the classification of subgroups of Suzuki groups 
 (\cite{HB}, Remark 3.12(e), p. 194) implies that each subgroup of $\Sz(q)$ 
 (except $\Sz(q)$ itself) 
has index $\ge q^2+1=\#(B)$. This implies that $\Sz(q)$ acts transitively on $B$. If the stabilizer $G_b$ of a point $b \in B$ has index $q^2+1$ then it follows easily from the same classification 
that $G_b$ is conjugate to the subgroup ${\mathfrak F}{\mathfrak H}$ generated by all $S(a,b)$ and $M(\lambda)$ and therefore the $\Sz(q)$-set $B$ is isomorphic to an {\sl ovoid} ${\mathcal O}=\Sz(q)/{\mathfrak F}{\mathfrak H}$ 
where the action of $\Sz(q)$ is known to be doubly transitive (\cite{HB}, Th. 3.3 on pp. 184--185 and
steps g) and i) of its proof on p. 187).
 By Remark \ref{St2}, this implies that the $\F_2[\Sz(q)]$-module $Q_B$ is absolutely simple. Recall that
$\dim_{\F_2}(Q_B)=\#(B)-1=q^2=2^{2(2k+1)}$.
By Theorem \ref{Sz}, there no absolutely simple nontrivial 
$\F_2[\Sz(q)]$-modules of dimension $< 2^{2(2k+1)}$. This implies that $Q_B$ is {\sl not} isomorphic to a tensor product of absolutely simple $\F_2[\Sz(q)]$-modules of dimension $>1$. Recall that all subgroups in $G=\Sz(q)$ (except $\Sz(q)$ itself) have index $\ge q^2+1> q^2=\dim_{\F_2}(Q_B)$. It follows from Corollary \ref{Very2} that the $G$-module $Q_B$ is very simple. Since $G \subset G'$, the $G'$-module $Q_B$ is also very simple. 
\end{proof}

\begin{proof}[Proof of Theorem \ref{redp}]
The cases (i) and (ii) of Theorem \ref{redp} follow 
 from Theorems \ref{L2QB} and \ref{SzQB} respectively 
applied to $G'=H$. 
\end{proof}

In light of Corollary \ref{mainvsc}  it would be interesting
to classify  all permutation subgroups $G \subset \Perm(B)$
with very simple $G$-modules $Q_B$.
We finish the paper by examples of very simple $Q_B$
 attached to Mathieu groups $\M_{11}$ and $\M_{12}$ and to 
 related group $\L_2(11)=\PSL_2(11)$.

\begin{thm}
\label{Atlas}
Let $n$ be a positive integer, $B$ a $n$-element set,
 $G \subset \Perm(B)$ a permutation group. Assume that $(n,G)$ enjoy one of the
following properties:
\begin{enumerate}
\item [(i)]
$n=11$ and $G$ is isomorphic either to $\L_2(11)$ or $\M_{11}$;

\item [(ii)]
$n=12$ and either $G \cong \M_{12}$ or
$G \cong \M_{11}$ and $G$ acts transitively on $B$. 
\end{enumerate}
Then the $G$-module $Q_B$ is very simple.
\end{thm}

\begin{proof}
Assume that $n=11$. Since $M_{11}$ contains a subgroup isomorphic to $\L_2(11)$  
(\cite{Atlas}, p. 18), it suffices to check the case of $G=\L_2(11)$,
in light  of Remark \ref{image}(iii). 

The group $G=\L_2(11)$  has two conjugacy classes of maximal subgroups 
of index $11$ and all other subgroups in $G$ have index greater than $11$ 
(\cite{Atlas}, p. 7). Therefore  
all  subgroups in $G$ (except $G$ itself)  have index greater than $10$ and the action of $G$ on the 
$11$-element set $B$ is transitive. 
 The permutation character (in both cases) is (in notations of \cite{Atlas}, p. 7) 
 $1+\chi_5$, i.e.,
$\bchi=\chi_5$. The restriction of $\chi_5$ to the 
set of $2$-regular elements coincides with absolutely irreducible
 Brauer character $\varphi_4$ (in notations of \cite{AtlasB}, p. 7). 
In particular, the corresponding $G$-module $Q_B$ 
 is absolutely simple and has dimension $10$.
Since $10=2 \cdot 5$ and $5$ is a prime, the very simplicity of 
the $G$-module $Q_B$ follows from Th. 5.4 of \cite{Zarhin}.
This proves the case (i).

Assume that $n=12$.

Suppose  $G=\M_{11}$ and the action of $G$ on the $12$-element 
 set $B$ is transitive.
All subgroups $G_b$ in $G$ of index $12$ are isomorphic to 
 $\L_2(11)$ (\cite{Atlas}, p. 18). 
It follows from the already proven case (i) for $\L_2(11)$ and  
 Remark \ref{oddeven}  that the $G$-module $Q_B$ is 
 very simple.

Suppose $G=\M_{12}$. The action of $G$ on the $12$-element $B$  is 
transitive, since all  subgroups in $G$ (except $G$ itself) have 
index $\ge 12$. All subgroups $G_b$ in $G$ of index $12$ are
isomorphic to $\M_{11}$ (\cite{Atlas}, p. 33). It follows from 
 the already proven case (i) for 
 $\M_{11}$ and Remark \ref{oddeven} 
 that the $G$-module $Q_B$ is very simple.
\end{proof}

Combining Corollary \ref{mainvsc} and Theorem \ref{Atlas} 
(with $B=\R, G=\Gal(f)$ and taking into account that the
 irreducibility of $f$ means that
 $\Gal(f)$ acts transitively on $\R$), 
 we obtain the following statement.

\begin{thm}
\label{mathieu}
Let $K$ be a field with $\fchar(K) \ne 2$,
 $K_a$ its algebraic closure,
$f(x) \in K[x]$ an irreducible separable polynomial of degree $n\ge 5$. Let $\R=\R_f \subset K_a$ be the set of roots of $f$, let $K(\R_f)=K(\R)$ be the splitting field of $f$ and $\Gal(f):=\Gal(K(\R)/K)$ the Galois group of $f$, viewed as a subgroup of $\Perm(\R)$.
Let $C_f$ be the hyperelliptic curve $y^2=f(x)$. Let  $J(C_f)$ be
its jacobian, $\End(J(C_f))$ the ring of $K_a$-endomorphisms of $J(C_f)$.

Assume that $n$ and $\Gal(f)$ enjoy one of the following properties:

\begin{enumerate}
\item[(i)]
$n=11$ and $\Gal(f)$ is 
isomorphic either to $\L_2(11)$ or to $\M_{11}$;
\item[(ii)]
$n=12$ and $\Gal(f)$ is isomorphic either to $\M_{11}$ or to 
$\M_{12}$;
\end{enumerate}

 Then either $\End(J(C_f))=\Z$ or 
$\fchar(K)>0$ and
$J(C_f)$ is a supersingular abelian variety.
\end{thm}

\end{document}